\documentclass[12pt, fleqn]{article}
\usepackage[cp1251]{inputenc}
\usepackage{latexsym,amsfonts,amssymb}
\usepackage{graphicx}

\usepackage{amsbsy}
\usepackage{amsmath}
\usepackage{epsf}

\newtheorem{theo}{Theorem}

\newcommand{\bt}{\begin{theo}}
\newcommand{\et}{\end{theo}}
\newcommand{\bd}{\begin{displaymath}}
\newcommand{\ed}{\end{displaymath}}

\newcommand{\be} {\begin{equation}}
\newcommand{\ee} {\end{equation}}
\newcommand{\ba}{\begin{array}{l}}
\newcommand{\ea} {\end{array}}
\newcommand{\bea}{\begin{eqnarray}}
\newcommand{\eea} {\end{eqnarray}}

\newcommand{\p} {\partial}

\begin{document}

\begin{center}
 {\Large \textbf{ Comments on the paper ``Solutions of Multitime Reaction-Diffusion PDE''} }
\medskip\\
{\bf Roman Cherniha~$^{\dag,\dag\dag}$ \footnote{\small
Corresponding author. E-mails: r.m.cherniha@gmail.com; roman.cherniha1@nottingham.ac.uk}}
 \\
{\it $^{\dag}$~Institute of Mathematics,  National Academy
of Sciences  of Ukraine,\\
 3, Tereshchenkivs'ka Street, Kyiv 01004, Ukraine \\
  $^{\dag\dag}$~School of Mathematical Sciences, University of Nottingham,\\
  University Park, Nottingham NG7 2RD, UK
}
 \end{center}

\begin{abstract}
The Comments are   devoted to the    paper ``Solutions of Multitime
Reaction-Diffusion PDE'' (Mathematics, {\bf 10} (2022), 3623), in
which main results are misleading and can be derived in a simple way
from those obtained earlier. Moreover, it is shown that  the exact
solution derived therein are simple generalizations of the known
solutions and are easily obtainable by the method of differential
constrains. The  Comments were  submitted  to the journal  Mathematics. 
\end{abstract}

Keywords: exact solution, diffusion-convection-reaction equation,
method of differential constrains,  scaling transformation.

\section{Multitime reaction-diffusion PDE} \label{sec-1}

The recent paper \cite{ghiu-ud-2022} is devoted to search for   exact solutions of a special class of PDEs, the so-called multitime reaction-diffusion equations

\be \label {1}
 h^1(t)\frac{\p u}{\p t^1}+\dots +  h^m(t)\frac{\p u}{\p t^m}=
 \mu \frac{\p^n u}{\p x^n} +  f\Big(u,\dots,\frac{\p^{n-1} u}{\p x^{n-1}}\Big),\ee
 where $u(t,x)$ is an unknown smooth function, $f$ and $ h^i(t)$  are  given  function from
  the class $C^1$, while  $\mu\not=0, \, t=(t^1,\dots,t^m)$.
   In the PDE theory, the class of PDEs   (\ref{1}) with  $n=2$ is known as
   the class of ultra-parabolic equations.

 From the very beginning, it should be stressed that the PDE class (\ref{1}) is not presented in
  a canonical form. As a result, the paper \cite{ghiu-ud-2022} consists of  many cumbersome
  statements and awkward formulae. Moreover, the technique used for constructing exact
   solutions is nothing else but a particular realization of
   {\it the  method of differential constrains}. This method was suggested by
   the know Soviet mathematician with the Ukrainian roots Yanenko
    in 1960s \cite{yanenko-64}.
   Although his pioneer work was never translated into English, one is well-known
    among experts in the field of exact solutions for nonlinear PDEs
    (see Chapter 5 in \cite{ch-se-pl-book} and references therein).
    Notably the method of compatible differential constrains (side conditions) \cite{olver-94}
    can be considered as
    a generalization of  the  method of differential constrains.

 Each PDE belonging to class (\ref{1}) is reducible to its canonical form
 \be \label {2}
 \frac{\p u}{\p \tau^1}+\dots+\frac{\p u}{\p \tau^m}= \mu \frac{\p^n u}{\p x^n} +  f\Big(u,\dots,\frac{\p^{n-1} u}{\p x^{n-1}}\Big).\ee
 The corresponding transformations is very simple and reads as 
 \be \label {3}
 \tau^i= \int\frac{dt^i}{ h^i(t)}, \  i=1,\dots,m
 \ee
 provided each function  $h^i(t)$ depends only on the variable $t^i$
 (equation (12)\cite{ghiu-ud-2022} contains exactly such functions) and   $h^i(t)\not=0$ (obviously this inequality  is true at least in some sub-domain where the functions $h^i(t)$ are continuously differentiable). If some  functions $h^i(t)$ depend on two or more variables then an appropriate transformation still exists and it is  shown in Appendix how that can be constructed.
 
 Transformation (\ref{3}) (or its analogue, see Appendix)  simultaneously reduce  equations (2), (3) and (4) from \cite{ghiu-ud-2022} to much simpler forms
 \be \label {2*}
 \frac{\p u}{\p \tau^1}+\dots+\frac{\p u}{\p \tau^m}= \mu \frac{\p^n u}{\p x^n} -k\frac{\p u}{\p x} +  f\Big(u,\dots,\frac{\p^{n-1} u}{\p x^{n-1}}\Big),\ee
  \be \label {4}
  \frac{\p w}{\p \tau^1}+\dots+\frac{\p w}{\p \tau^m}= 0 \ee
  and
  \be \label {5}
  \frac{\p v}{\p \tau^1}+\dots+\frac{\p v}{\p \tau^m}= 1, \ee
  respectively.
  According to  the method of  differential constrains,
  the {\it linear} differential equations (\ref{4}) and (\ref{5})  are additional constrains
  used for constructing exact solutions of {\it the nonlinear PDE} (\ref{1}).
  Notably,  the idea of application of linear differential equations for solving
  nonlinear PDEs was used earlier in many papers. For example,
  the method of additional generating conditions \cite{ch-96,ch-98} is based on additional
   constrains
  in the form of linear ordinary differential equations.

 Thus, one should study the class of PDEs (\ref{2}) with the additional constrains
 (\ref{4}) and (\ref{5})  instead of  PDEs used in  \cite{ghiu-ud-2022}
 (see formulae (1), (3) and (4) therein).

 Now it can be easily demonstrated that Propositions 1--4, which form the theoretical core of  the technique developed in \cite{ghiu-ud-2022}, are rather trivial statements.
 Consider, for example, Proposition 1, which simply says that the function
 \be \label {6}
 u=\phi\Big( \tau, \,x+kv(\tau)+P_1(\tau) \Big)\ee
 is an exact solution of Eq.(\ref{2}) provided the functions $\phi(\tau,x)$,  $P_1(\tau)$
 and $v(\tau)$  are the exact solutions of Eqs.(\ref{2*}), (\ref{4}) and  (\ref{5}),
  respectively. Of course, sufficient smoothness of the above functions in the relevant
  domains is assumed as well.
 On the other hand, if one substitutes formula (\ref{6}) into  Eq.(\ref{2})
 and calculates derivatives then immediately  Eq.(\ref{2*}) for the function $\phi$
 is obtained. But the latter is a solution of  Eq.(\ref{2*}) by the above assumption.
 As a result, the authors propose to find exact solutions of Eq.(\ref{2}) using the given
  solution of  Eqs.(\ref{2*}), which is more complicated. It is a useless technique.

 However, using   the  method of differential constrains, the more constructive  result
  is obtainable in a straightforward  way. In fact, Eq.(\ref{2})  can be rewritten
   in the form
 \be \label {7}
 \frac{\p u}{\p \tau^1}+\dots+\frac{\p u}{\p \tau^m}= \mu \frac{\p^n u}{\p y^n} -k\frac{\p u}{\p y} +f\Big(u,\dots,\frac{\p^{n-1} u}{\p y^{n-1}}\Big)\ee
 by introducing the variable $y=x+k\tau^m$. Now we consider the linear PDE
  \be \label {8}
  \frac{\p u}{\p \tau^1}+\dots+\frac{\p u}{\p \tau^m}= 0 \ee
  as the differential constrain for Eq.(\ref{2*}).
  The general solution of  Eq.(\ref{8}) is readily  constructed:
  \be \label {9}
 u=\psi(\omega_1,\dots,\omega_{m-1},y) \ee
 where $\psi$ is an arbitrary smooth function,  $\omega_j=\tau^m -\tau^j, \, (j=1,\dots,m-1)$.

 Thus, formula (\ref{9})  is an exact solution of Eq.(\ref{2*}) provided the function
  $\phi(\tau, y)$ (with the  variables $\tau=(\tau^1,\dots,\tau^m)$
  to-be-determined as parameters) is a  solution of
  the {\it ordinary differential equation} (ODE)
  \be \label {10}
 \mu \frac{d^n u}{d y^n} -k\frac{d u}{d y} +f\Big(u,\dots,\frac{d^{n-1} u}{d y^{n-1}}\Big) = 0.\ee
 Although (\ref{10}) is still the nonlinear ODE,  many such equations
 with correctly-specified function $f$ are integrable  and the corresponding exact solutions
 can be found either in the well-known handbooks like \cite{kamke},
 or by using computer algebra solvers from  Maple, Mathematica etc.
 Finally, each exact solution  (\ref{9}) of  Eq.(\ref{2*}) is automatically the solution
  \be \label {11}
 u=\psi(\omega_1,\dots,\omega_{m-1},x+k\tau^m) \ee
 of the nonlinear PDE (\ref{2}). In the case $m=1$, one obtains the plane wave solutions,
 which are common in real-world  applications provided they are bounded and
 non-negative, i.e., they are traveling fronts \cite{gild-ker-04,ch-se-pl-book}.

 Proposition 2 in \cite{ghiu-ud-2022} is a straightforward  generalization of Proposition 1.
  On the other hand, there is no need to introduce $r+1$ particular solutions
  of Eg.(\ref{4}) in order to generalize Proposition 1  because  Eg.(\ref{4}) is
   integrable and its general solution can be simply used.
   Moreover, a  $r$-parameter family of exact solutions of Eg.(\ref{2*}) is needed
    in order to obtain the relevant family of solutions of Eg.(\ref{2}).
     Obviously, this statement cannot be used for real applications because
      the authors again  try to solve the  given PDE using solutions of the  more
      complicated PDE. On the other hand, using the general solutions of Egs.(\ref{4})
 and (\ref{5}),  formula (\ref{6})  can be rewritten in the form
 \be \label {12}
 u=\phi\Big( \tau, \,x+k\tau^m+P(\omega_1,\dots,\omega_{m-1}) \Big)\ee
 where $P(\omega_1,\dots,\omega_{m-1})$ is an arbitrary smooth function and
 $\phi( \tau, \,x)$ is a solution of Eg.(\ref{2*}). In particular,  one obtains the plane wave
 solutions in the case $m=1$ provided   $\phi( \tau, \,x)$ does not depend on $\tau=\tau^1$.

 \section{Exact solutions of a specific multitime reaction-diffusion equation } \label{sec-2}

Exact solutions in explicit forms are presented  in Section 4 \cite{ghiu-ud-2022}  only for the nonlinear PDE
\be \label {2-1}
 t^1\frac{\p u}{\p t^1}+\dots +  t^m\frac{\p u}{\p t^m}= \mu \frac{\p^2 u}{\p x^2} -au^3 + bu^2 \ee
 with the constants $\mu >0, \, a>0$ and $b\not=0$.
 All the exact solutions of  PDE  (\ref{2-1}) constructed in \cite{ghiu-ud-2022} are nothing else but a direct generalizations of the known solutions of this equation in the case $m=1$. Moreover, the solutions  in the case $m>1$ are obtainable from those with $m=1$ in a straightforward way using the method of differential constrains.

 First of all,  PDE  (\ref{2-1}) is reducible to its canonical form
 \be \label {2-2}
 \frac{\p u}{\p \tau^1}+\dots+\frac{\p u}{\p \tau^m}= \mu \frac{\p^2 u}{\p x^2}  -au^3 + bu^2 \ee
 by the transformation
 \be \label {2-3}
 \tau^i= \ln t^i, \  i=1,\dots,m.
 \ee
 Further, it can be noted that PDE  (\ref{2-2}) can be simplified to the form
 \be \label {2-4}
 \frac{\p u_*}{\p \tau_*^1}+\dots+\frac{\p u_*}{\p \tau_*^m}=
 \frac{\p^2 u_*}{\p x_*^2}  -u_*^3 + u_*^2 \ee
 using the scaling transformations
  \be \label {2-5}
  u_*=\frac{a}{b}u, \, x_*=\frac{|b|}{\sqrt{\mu a}}x, \, \tau_*^i=\frac{b^2}{a}\tau^i \quad  (i=1,\dots,m).\ee
  So, one can examine the  nonlinear PDE (the symbol star is omitted in what follows)
   \be \label {2-2*}
 \frac{\p u}{\p \tau^1}+\dots+\frac{\p u}{\p \tau^m}= \frac{\p^2 u}{\p x^2}  -u^3 + u^2 \ee
 instead of (\ref{2-1}) without losing a generality.

 Now one realizes that PDE (\ref{2-2*}) with $m=1$  is the known Huxley equation  and
 its exact solutions  were  constructed, for example,   in \cite{cla-mans-94, a-h-b-94}
 using non-classical symmetries.
 This equation can be also thought as  a particular case  of
 the famous  Fitzhugh-Nagumo equation
 \be \label {2-6} \frac{\p u}{\p \tau}= \frac{\p^2 u}{\p x^2} +  u(1-u)(u-\delta). \ee
 Exact solutions of the latter were constructed in many studies
 \cite{kawa-tana-83, cla-mans-94, a-h-b-94}. In the book
 \cite{gild-ker-04}, traveling waves of Eqs.(\ref{2-6}) are
 presented.
 It can be noted that all the exact solutions constructed  in Section 4 \cite{ghiu-ud-2022}
  is nothing else but direct generalizations of those of the Huxley equation.

 Let us show this. First of all, it should be noted that an arbitrary solution $u(\tau,x)$
  of PDE (\ref{2-2*}) can be multiplied to the form
 $u(\tau+\tau_0,\pm x+x_0)$ with an arbitrary constant vector $\tau_0$ and a constant $x_0$.
  It is a trivial consequence of invariance of   PDE (\ref{2-2*}) w.r.t. the time and space
  translations and the discrete transformation $x \to -x$.  In particular it means that
   exact solutions $u_1, \, u_2$ and $u_3$  are equivalent to those
    $u_2, \, u_4$ and $u_6$, respectively  (see pages 9--10 in \cite{ghiu-ud-2022}).
     Notably the authors should present also the solution $u_8$ replacing
      $x$ by $ -x$ in $u_7$.

 The most interesting from the   applicability point of view are the solutions
 $u_1$ and $u_3$ \cite{ghiu-ud-2022},  which   are reducible to much simpler form (see the above transformations):
 \be \label {2-7}
 u= \frac{\sqrt{2}}{x-\sqrt{2}\tau^m + P(\omega_1,\dots,\omega_{m-1}) } \ee
 and
  \be \label {2-8}
  u= \frac{1}{1+ P(\omega_1,\dots,\omega_{m-1})\exp\Big(\frac{x}{\sqrt 2}-\frac{\tau^m}{2}\Big)}, \ee
respectively. Here  $P(\omega_1,\dots,\omega_{m-1})$ is an arbitrary
smooth function of the variables
  $\omega_j=\tau^m -\tau^j, \, (j=1,\dots,m-1)$.

   In the case $m=1$,  solution  (\ref{2-7}) takes the form
   \be \label {2-9}
 u=\frac{\sqrt{2}}{x-\sqrt{2}\tau +x_0} , \ee
 which is nothing else but a very particular case of the exact solution
 (5.25)\cite{cla-mans-94}. Similarly,  solution  (\ref{2-8}) reduces to
 \be \label {2-10}
  u= \frac{1}{1+ C\exp\Big(\frac{x}{\sqrt 2}-\frac{\tau}{2}\Big)}. \ee
 Assuming that the constant $C= e^{x_0}>0$, one arrives at the   traveling front
  \be \label {2-11}
  u= \frac{1}{2}{1+ \tanh\Big(\frac{-\sqrt{2}}{4}(x+x_0)+\frac{1}{4}\tau\Big)}, \ee
 which was identified in \cite{kawa-tana-83}. In the case $C<0$, the
 exact solution (\ref{2-11}) with the function $\coth$ is obtained.

 Finally, it should be stressed that  the method of differential constrains with using
 the linear PDE (\ref{8}) as the given constrain allow us easily  to
 identify a set of exact solutions of  PDE (\ref{2-2*}), which are
 reducible to those obtained in \cite{ghiu-ud-2022}. Actually, the
 authors   do not apply Propositions 1--4 for finding exact solutions
 but ODE (15)--(16) \cite{ghiu-ud-2022} are used, therefore they  implicitly
 use  the method of differential constrains for finding exact
 solutions. By the way, ODE (16) \cite{ghiu-ud-2022} is fully
 integrable in terms of elliptic functions, hence its  general
 solution involving  $y_5, \, y_6$ and $y_7$ \cite{ghiu-ud-2022} as particular cases
  could be presented therein.

  In conclusion, I would like to stress that nowadays  there are many papers
  devoted
  to search for exact solutions of nonlinear PDEs, in which new
  methods are suggested without knowing the state-of-art. As a result,
  the newly suggested
  methods very often are not new and the solutions obtained are 
  straightforward generalizations of those derived earlier. Study
  \cite{ghiu-ud-2022} is a typical example.

 \section{Appendix}
 If  PDE (\ref{1})  involves the smooth  functions $ h^i(t), \ i=1,\dots,m$ of the more general form than it was assumed above then a transformation reducing  PDE (\ref{1}) to PDE (\ref{2})  still exists. Let us construct the transformation in the case $m=2$ in order to avoid cumbersome formulae (the case $m>2$ can be examined in the same way). Assuming that a non-degenerate  transformation in question  possesses the form 
 \be \label {A-1}  \tau^1=H^1(t^1,t^2),  \, \tau^2=H^2(t^1,t^2), \ee
 one easily calculates 
 \be \label {A-2}  h^1(t)\frac{\p u}{\p t^1}+  h^2(t)\frac{\p u}{\p t^2}=
 \Big(h^1(t)\frac{\p H^1}{\p t^1} + h^2(t)\frac{\p H^1}{\p t^2} \Big)\frac{\p u}{\p \tau^1} +
  \Big(h^1(t)\frac{\p H^2}{\p t^1} + h^2(t)\frac{\p H^2}{\p t^2} \Big)\frac{\p u}{\p
  \tau^2}.  
 \ee
 So, the local transformation (\ref{A-1})  reduces PDE (\ref{1}) to PDE (\ref{2}) with $m=2$ provided the functions $H^1$ and $H^2$ form an arbitrary solution of the linear decoupled system of the first-order PDE
 \be  \label {A-3} \ba \medskip
 h^1(t)\frac{\p H^1}{\p t^1} + h^2(t)\frac{\p H^1}{\p t^2}=1, \\
 h^1(t)\frac{\p H^2}{\p t^1} + h^2(t)\frac{\p H^2}{\p t^2}=1.
\ea \ee
 According to the theory of the linear first-order PDEs the solutions of (\ref{A-3}) possess the form 
 \be \label {A-4}\ba \medskip
H^1= \int\frac{dt^1}{ h^1(t)} + W^1(I), \\
H^2= \int\frac{dt^1}{ h^1(t)} + W^2(I),
\ea \ee
where $ W^1$ and $ W^2$  are functions of the   first integral $I(t^1,t^2)$ of the ODE 
\be \label {A-5}  \frac{dt^1}{h^1(t)}=  \frac{dt^2}{h^2(t)}.
\ee
In the case of real-world applications, the functions $h^1(t)$ and $h^2(t)$ are such that  the first integrals of ODE (\ref{A-4}) can be written down in an explicit form, otherwise those  can be presented in an implicit form. Obviously, the functions  
$ W^1$ and $ W^2$ should be linearly independent, otherwise a degenerate transformation is obtained.

\end{document}